# ON THE POISSON EQUATION AND DIFFUSION APPROXIMATION 3

By E. Pardoux and A. Yu. Veretennikov

*Université de Provence and University of Leeds*

We study the Poisson equation $Lu + f = 0$ in $\mathbb{R}^d$, where $L$ is the infinitesimal generator of a diffusion process. In this paper, we allow the second-order part of the generator $L$ to be degenerate, provided a local condition of Doeblin type is satisfied, so that, if we also assume a condition on the drift which implies recurrence, the diffusion process is ergodic. The equation is understood in a weak sense. Our results are then applied to diffusion approximation.

**1. Introduction.** This is the third in a series of papers devoted to the study of the Poisson equation in $\mathbb{R}^d$ and diffusion approximation. In this paper we consider the degenerate case.

The study of diffusion approximation [i.e., obtaining the limit of $Y^\varepsilon$ in (19)] was initiated by Khasminskii [5], and developed by many authors, including Papanicolaou, Stroock and Varadhan [10] and Kushner [8]. Such results, and the formulation of the limiting stochastic differential equation, require the solution of a Poisson equation $Lu + f = 0$, where $L$ is the infinitesimal generator of a Markov process (at least in the case where the disturbance is Markovian; in the non-Markov case a substitute of the Poisson equation replaces it), whose right-hand side $f$ is the highly oscillating coefficient of the approximating differential system. When the disturbance in the approximating ODE is compact valued, the Poisson equation is formulated in a compact set, and the corresponding theory is well known; the result can be proved under quite explicit conditions on the coefficients (see [2], Chapter 12, Section 2). When, however, the disturbance of the diffusion takes values in all of $\mathbb{R}^d$, there was until recently no way of deriving estimates for the solution of the Poisson equation in terms of explicit conditions on the data; see Chapter 12, Section 3 in [2].

This was the starting point for our work. We focused on the case where the disturbance is an ergodic $\mathbb{R}^d$-valued diffusion process, and the ergodicity follows from explicit conditions on the coefficients. In the first paper [13], using mainly probabilistic arguments (together with some estimates from the theory of partial differential equations), we solved the Poisson equation in $\mathbb{R}^d$ for the generator of an elliptic and ergodic diffusion, and obtained estimates (which we believe are rather sharp) of the solution. We then used that result to establish a diffusion approximation result under very explicit conditions on the coefficients. In the second paper [14], we considered the case where the coefficients of the equation (18) depend on the $Y^\varepsilon$ process. This forced us to study a Poisson equation where both the PDE operator and the right-hand side depend on a parameter, and establish regularity results of the solution in terms of that parameter. We were forced for that purpose to use essentially results from the PDE theory.

The aim of this third paper is to consider the situation of the first paper, where we now relax the ellipticity assumption. While a condition on the behavior of the drift at infinity [condition $(A_b)$] implies the positive recurrence, irreducibility, which was in our previous works a consequence of ellipticity, is now a consequence from a type of "local Doeblin condition" [condition $(D_\ell)$]. While those conditions are not explicit conditions on the coefficients of the diffusion, they are implied both by the ellipticity assumption and by the "restricted Hörmander condition" (i.e., the assumptions that the diffusion vector fields, together with their brackets of arbitrary order, span the whole space at each point). We further give one example where none of these conditions holds, while our condition $(D_\ell)$ is satisfied, together with the additional "regularity" condition $(A_T)$.

We then consider a weak formulation of the Poisson equation, which is solved by the same probabilistic formula as in the elliptic case. We finally apply those results to the diffusion approximation problem. We prove weak convergence in the sense of the $S$-topology of Jakubowski [4]. The difficulty in proving convergence in a stronger sense is related to the lack of smoothness of the solution of the Poisson equation. Let us also mention that our diffusion approximation is "less general" than the one considered in [13] (except that we relax the ellipticity condition, as explained above), in that the approximating differential system does not contain a stochastic integral (i.e., the coefficient $H$ in [13] does not appear here).

Let us point out that precise regularity of the solution of the Poisson equation under the Hörmander condition follows from Theorem 18 in [16].

The paper is organized as follows. Section 2 contains our assumptions, as well as some essential ergodicity results from Veretennikov [23]. Section 3 is devoted to the study of one example of degenerate coefficients, which satisfy our assumptions. The Poisson equation is studied in Section 4, while the diffusion approximation result is derived in Section 5.

**2. Moment bounds and convergence to the invariant measure.** Consider the stochastic Itô equation

$$dX_t = b(X_t)\,dt + \sigma(X_t)\,dB_t, \qquad X_0 = x \in \mathbb{R}^d, \tag{1}$$

where $\{B_t, t \geq 0\}$ is a $k$-dimensional Brownian motion, $b$ is a locally Lipschitz vector-function of dimension $d$ and $\sigma$ is a $d \times k$ matrix-valued locally Lipschitz function. We assume that $\sigma\sigma^*$ is bounded and possibly degenerate, and that the unique solution of (1) satisfies

$$(A_T) \qquad \forall R > 0 \qquad \sup_{|x| \leq R} \mathbb{E}_x \inf\{t \geq 0 : |X_t| \geq R+1\} < \infty.$$

Let us introduce the following recurrence condition:

$$(A_b) \qquad \lim_{|x| \to \infty} (b(x), x) = -\infty.$$

Note that this condition prevents the solution of the SDE (1) from exploding, so that the process $\{X_t\}$ is well defined for all $t > 0$. Let $R > 0$ and $\tau = \tau^R = \inf(t \geq 0 : |X_t| \leq R)$.

Finally, we assume the following "local Doeblin" type condition.

Let $B \subset \mathbb{R}^d$ and

$$\tau_0^B := \inf\{t \geq 0 : X_t \in B\},$$

and (in the following formula, $t_B > 0$ depends only on $B$)

$$\tau_{\ell+1}^B := \inf\{t \geq \tau_\ell^B + t_B : X_t \in B\}.$$

Define the "*process in $B$*" in discrete time as $X_n^B := X_{\tau_n^B}$. Denote by $P^B(n, x, dx')$ the $n$-step transition probability of $(X_n^B)$. We say that the *local Doeblin condition* holds true for the process $\{X_t\}$ if for any $R' > 0$ there exists $R > R'$ such that the process in $B = B_R := \{x \in \mathbb{R}^d : |x| \leq R\}$ satisfies the following: there exists an integer $n_0 = n_0(R) > 0$ such that

$$(D_\ell) \quad \inf_{|x|,|x'| \leq R} \int_B \min\left\{\frac{P^B(n_0, x, dx'')}{P^B(n_0, x', dx'')}, 1\right\} P^B(n_0, x', dx'') =: q(R, n_0) > 0,$$

where $\frac{P^B(n_0, x, dx'')}{P^B(n_0, x', dx'')}$ is defined as follows. Let

$$P^B(n_0, x, dx'') = \varphi_{x,x'}(x'') P^B(n_0, x', dx'') + \nu_{x,x'}(dx'')$$

be the decomposition of $P^B(n_0, x, dx'')$ into its absolutely continuous part w.r.t. $P^B(n_0, x', dx'')$, and the part $\nu_{x,x'}(dx'')$ which is singular w.r.t. $P^B(n_0, x', dx'')$. Then

$$\frac{P^B(n_0, x, dx'')}{P^B(n_0, x', dx'')} := \varphi_{x,x'}(x'').$$

The assumption $(D_\ell)$ requires, in particular, that the mass of the singular part is not close to 1, and moreover, it imposes a certain quantitative estimate on the total variation norm for the difference of two measures uniformly on the compact $B$. We assume throughout the paper that the process $\{X_t\}$ satisfies this assumption. We shall give in the next section one example with a nonelliptic diffusion coefficient, for which conditions $(A_T)$ and $(D_\ell)$ hold. The proof that example satisfies our conditions will use stronger conditions, which are easier to verify.

We note that the two assumptions $(A_b)$ and $(D_\ell)$ imply the existence and uniqueness of an invariant probability measure. For this and the proof of the next proposition, see [23]. Note that in [23] the "process in $B$" is defined in a slightly different manner, since it is extracted from the sequence $\{X_n, n \in \mathbb{N}\}$ rather than from $\{X_t, t \geq 0\}$ as defined here. However, the adaptation of those proofs is rather obvious.

PROPOSITION 1. *Under the assumptions $(A_T)$, $(A_b)$ and $(D_\ell)$, for all $m' > m + 2 > 2$, there exists $C$ such that for all $x \in \mathbb{R}^d$, $t > 0$,*

$$\mathbb{E}_x |X_t|^m \leq C(1 + |x|^{m'}). \tag{2}$$

*Moreover,*

$$\mathbb{E}_\mu |X_t|^m < \infty \qquad \forall\, m > 0, \tag{3}$$

*and for any $k > 0$, $2k + 2 < m$,*

$$\mathrm{var}(\mu_t^x - \mu) \leq C(1 + |x|^m)(1 + t)^{-(k+1)}, \tag{4}$$

*where "var" denotes the total variation norm of a signed measure over the Borel $\sigma$-field, $\mu_t^x$ is the law of $X_t$ when $X_0 = x$, $\mu$ is the unique invariant measure of $X$ and $\mathbb{E}_\mu$ means the expectation w.r.t. $\mu$.*

PROPOSITION 2. *Let the assumptions $(A_T)$, $(A_b)$ and $(D_\ell)$ be satisfied. Then for any $p > 0$,*

$$\mathbb{E}_x \bigg( \sup_{0 \leq t' \leq t} |X_{t'}|^p \bigg) = o(\sqrt{t}) \qquad as\ t \to \infty.$$

The proof of Proposition 2 is similar to that in [13]; hence, we drop it. The following corollary will be used in Section 4, for the proof of tightness.

COROLLARY 1. *Under the same assumptions, for any $T > 0$, $p > 0$,*

$$\varepsilon \mathbb{E}_x \bigg( \sup_{0 \leq t \leq T} |X_{t/\varepsilon^2}|^p \bigg) \to 0 \qquad as\ \varepsilon \to 0.$$

**3. Sufficient conditions and one example.** In this section, we first state two conditions, which we prove to be, respectively, stronger than $(A_T)$, and stronger than $(D_\ell)$. Then we give one example with a degenerate diffusion coefficient, which satisfies those stronger conditions.

PROPOSITION 3. *A sufficient condition for condition $(A_T)$ to hold is that for each $R > 0$, there exists $f \in C(\mathbb{R}^d, [0,1])$ with $\mathrm{supp}(f) \subset \{x; |x| \geq R+1\}$ and $t_0$ such that*

$$\inf_{|x| \leq R} \mathbb{E}_x f(X_{t_0}) > 0.$$

PROOF. Let $R > 0$. It follows from our condition that there exists $c, t_0 > 0$ such that for all $x \in \mathbb{R}^d$, $|x| \leq R$, $P_x(|X_{t_0}| \geq R+1) \geq c$; this implies $P_x(|X_{\{t_0\}}| \geq R|) \geq c$, and the same is true for $R+1$ instead of $R$, with new $c$ and $t_0$. Let

$$S_{R+1} = \inf\{t; |X_t| \geq R+1\}.$$

It follows from the Markov property of $\{X_t, t \geq 0\}$ and the previous estimate that for all $x \in \mathbb{R}^d$,

$$\mathbb{P}_x(S_{R+1} > nt_0) \leq (1-c)^n.$$

The result follows. $\square$

We now formulate what we call the condition $(D_{s\ell})$ ("strong local Doeblin condition").

For each $R > 0$, there exists $A_R \subset B_R$, $t_R > 0$, $c(R) \geq 1$ such that for all $x \in B_R$, the transition probability of our diffusion process $\{X_t, t \geq 0\}$ satisfies

$$p(t_R, x; dy) = q(t_R, x, y)\mu(dy) + \nu(t_R, x; dy),$$

$(D_{s\ell})$

$$\frac{1}{c(R)} \leq q(t, x, y) \leq c(R), \qquad y \in A_R,$$

where $\mu$ is a probability measure on $\mathbb{R}^d$ such that $\mu(A_R) > 0$.

Note that the upper bound here is not actually necessary for our aims; however, it holds along the lower bound in all cases known to the authors.

PROPOSITION 4. *The strong local Doeblin condition $(D_{s\ell})$ implies the local Doeblin condition $(D_\ell)$.*

PROOF. We choose an arbitrary $R > 0$, and denote $B = B_R$. We decompose the transition probability of the process in $B$ (defined with $t_B = t_R$) as follows. For $x \in B$,

$$P^B(1, x, dx') = \mathbb{P}_x(X_{t_R} \in dx', X_{t_R} \in B) + \nu'(x, dx').$$

It follows from $(D_{s\ell})$ that

$$\mathbb{P}_x(X_{t_R} \in dx', X_{t_R} \in B) \geq \frac{1}{c(R)} \mathbb{1}_{A_R}(x')\mu(dx').$$

Hence $(D_\ell)$ holds with $n_0 = 1$ and $q(R, 1) = \frac{\mu(A_R)}{c(R)}$. □

EXAMPLE 1. Let $b$ and $\sigma_0$ satisfy the above conditions $(A_T)$ and $(A_b)$, $b \in C^1(\mathbb{R}^d, \mathbb{R}^d)$, $\sigma_0 \in C^2(\mathbb{R}^d, \mathbb{R}^{d \times d})$, and let $\sigma_0$ be uniformly nondegenerate. Let $\alpha : \mathbb{R}^d \to [0, 1]$ be a $C^1$ mapping, such that the set $\{\alpha = 0\}$ is the union of countably many disjoint connected closed subsets of $\mathbb{R}^d$, such that each bounded subset of $\mathbb{R}^d$ intersects at most finitely many of those, and the set $\{\alpha > 0\}$ is connected. We now assume that for some $\delta > 0$ such that $\{\alpha > \delta\}$ is connected and for each $R' > 0$ there exists $R > R'$ such that the set $\{|x| = R\}$ does not intersect the set $\{\alpha \leq \delta\}$, and moreover that for any $R > 0$, there exists $M$ such that the solution of

$$\frac{dx}{dt}(t) = b(x(t))$$

exits in time less than $M$ from $\{\alpha \leq \delta\}$, whenever $x(0) \in \{\alpha < \delta\} \cap B_R$. Let $\sigma(x) = \alpha(x)\sigma_0(x)$. Then the pair $(b, \sigma)$ satisfies the assumptions $(A_T)$ and $(D_\ell)$.

We first prove the following.

LEMMA 1. *The condition of Proposition 3 is satisfied in Example 1.*

PROOF. We consider the stochastic equation for the process $\{X_t, t \geq 0\}$, written in Stratonovich form (the reason for this is that we shall soon use Stroock and Varadhan's support theorem), that is,

$$dX_t = \tilde{b}(X_t) \, dt + \sigma(X_t) \circ dB_t,$$

where

$$\tilde{b}_i(x) = b_i(x) - \frac{1}{2}\left(\sum_{j,k} \frac{\partial \sigma_{ij}}{\partial x_k} \sigma_{kj}\right)(x).$$

We will use below the notation $(\nabla \sigma)\sigma_0$ for the vector $2\alpha^{-1}(b - \tilde{b})$. We now consider the controlled ODE

$$\frac{dy}{dt}(t) = \tilde{b}(y(t)) + \sigma(y(t))u(t),$$
$$y(0) = X_0,$$

where we choose the feedback control $u(t) = \Phi(y(t))$, with

$$\Phi(x) = \begin{cases} \frac{1}{2}\sigma_0^{-1}(\nabla \sigma)\sigma_0(x), & \text{if } \alpha(x) > 0, \\ 0, & \text{if } \alpha(x) = 0. \end{cases}$$

It is easy to check that $\{y(t), t \geq 0\}$ coincides with the solution of the ODE

$$\frac{dx}{dt}(t) = b(x(t)), \qquad x(0) = X(0).$$

Let

$$\tau = \tau(x(0)) = \inf\{t > 0; \alpha(x(t)) \geq \delta\}.$$

Choose $\rho = (2\|\nabla\alpha\|_{\infty,R})^{-1}\delta$, $g \in C(\mathbb{R}^d; [0,1])$ with $\mathrm{supp}(g) \subset \{x, |x| \leq \rho\}$ and $g \equiv 1$ on a neighborhood of 0. We have used the notation

$$\|\nabla\alpha\|_{\infty,R} = \sup_{|x| \leq R} \|\nabla\alpha(x)\|.$$

The above considerations and Stroock and Varadhan's support theorem (cf. [17]) imply that

$$\mathbb{E}_{x(0)} g(X_\tau - x(\tau)) > 0.$$

Moreover, that last quantity depends continuously on $x(0)$, hence it is bounded away from zero for $x(0) \in \{\alpha \leq \delta\} \cap B_R$. Hence our construction yields that with a probability which is bounded away from zero, $\alpha(X_\tau) \geq \delta/2$, where $\tau \leq M$ is a deterministic time which depends only on $x(0)$. $\tau = 0$ whenever $\alpha(x(0)) \geq \delta$.

Let $f \in C(\mathbb{R}^d, [0,1])$ satisfy $\mathrm{supp}(f) \subset \{x; |x| \geq R+1\}$ and $f(x) = 1$, whenever $|x| \geq R + 2$. Using again Stroock and Varadhan's support theorem, we have that

$$\inf_{x \in B_R} \inf_{x' \in B_R \cap \{\alpha \geq \delta/2\}} \mathbb{E}_{x'} f(X_{2M-\tau(x)}) > 0.$$

Proposition 3 with $t_0 = 2M$ and the above $f$ now follows from the Markov property. □

LEMMA 2. *The pair $(b, \sigma)$ from Example 1 satisfies the condition $(D_{s\ell})$.*

PROOF. It follows from the proof of Lemma 1 that there exists $\xi > 0$ and a mapping $\tau \in C(B_R, [0, M])$ such that for all $x \in B_R$,

$$\mathbb{P}_x(\alpha(X_{\tau(x)}) \geq \delta/2) \geq \xi.$$

Next, we choose a closed ball $A \subset \mathrm{int}\, B_R \cap \{\alpha > 0\}$. Using again Stroock and Varadhan's support theorem, we deduce that there exists $N > M$ such that

$$\inf_{x \in B_R} \inf_{x' \in B_R \cap \{\alpha \geq \delta/2\}} \mathbb{P}_{x'}(X_{N-\tau(x)} \in A) > 0.$$

Combining the above two statements with the help of the Markov property, we obtain that

(5) $$\inf_{x \in B_R} \mathbb{P}_x(X_N \in A) > 0.$$

Next we choose another closed ball $A'$ such that $A \subset \text{int}\, A' \subset A' \subset B_R \cap \{\alpha > 0\}$.

For any function $\varphi \in C(A', \mathbb{R}_+)$, with $\text{supp}(\varphi) \subset A$, we consider the solution $\{u(t,x), 0 \leq t \leq 1, x \in A'\}$ of the backward linear parabolic PDE [here $a(x) = \sigma \sigma^*(x)$]

$$\frac{\partial u}{\partial t}(t,x) + \frac{1}{2}\sum_{ij} a_{ij}(x)\frac{\partial^2 u}{\partial x_i \partial x_j}(t,x) + \sum_i b_i(x)\frac{\partial u}{\partial x_i}(t,x) = 0,$$

$$0 < t < 1, x \in A',$$

$$u(1,x) = \varphi(x), \qquad u(t,x) = 0, \qquad x \in \partial A'.$$

We have that for all $x \in A'$, $u(0,x) = \mathbb{E}_x \varphi(Y_1)$, where the process $\{Y_t, 0 \leq t \leq 1\}$ is the solution of the SDE (1), which is killed when it reaches the boundary of the set $A'$. It follows from the parabolic Harnack inequality (see, e.g., [7], page 131) that there exists $N > 0$ such that

$$\sup_{x,x' \in A} \frac{u(0,x)}{u(0,x')} \leq N,$$

that is,

$$\sup_{x,x' \in A} \frac{E_x \varphi(Y_1)}{E_{x'}\varphi(Y_1)} \leq N,$$

for all $\varphi \in C(A', \mathbb{R}_+)$, with $\text{supp}(\varphi) \subset A$. We choose one particular point $x_0 \in A$, and define $\mu(dy) = \mathbb{P}_{x_0}(Y_1 \in dy)$. It follows from the above that for each $x \in A$, $\mathbb{P}_x(Y_1 \in dy)$ is absolutely continuous with respect to $\mu$, and moreover the Radon–Nikodym derivative $q(x,y)$ satisfies

(6) $$N^{-1} \leq q(x,y) \leq N,$$

for all $x, y \in A$.

Condition ($D_{s\ell}$) now follows from (5), (6) and the Markov property. □

REMARK 1. It is rather clear that one can verify our assumptions in many other situations, where $\det(a(x))$ may vanish in a similar fashion as $\alpha(x)$ does in Example 1. All that is to be verified is a condition like (5), both for a set of the same type as $A$ and for $B^c_{R+1}$.

REMARK 2. In the strictly elliptic case the same arguments based on Harnack's inequality establish the condition ($D_{s\ell}$), provided $a = \sigma\sigma^*/2$ is continuous, and $b$ locally bounded. The same is true, with $\mu = $ Lebesgue measure, whenever the coefficients are smooth, and the Lie algebra of vectors fields generated by the columns of the matrix $\sigma$ has full rank at any point of $\mathbb{R}^d$.

**4. The Poisson equation in $\mathbb{R}^d$.** We consider the Poisson equation in $\mathbb{R}^d$

$$Lu(x) = -f(x), \tag{7}$$

where

$$L = \sum a_{ij}(x)\partial_{x_i}\partial_{x_j} + \sum b_i(x)\partial_{x_i},$$

with

$$a(x) = \sigma\sigma^*(x)/2,$$

and $f \in C(\mathbb{R}^d)$ satisfies

$$|f(x)| \leq C(1+|x|)^\beta \qquad \text{for some } \beta \in \mathbb{R},$$

so that due to Proposition 1, $f$ is integrable with respect to the invariant measure $\mu$, and

$$\int f(x)\mu(dx) = 0. \tag{$A_f$}$$

In the nondegenerate case, the solution of (7) has the representation

$$u(x) = \int_0^\infty \mathbb{E}_x f(X_s)\,ds. \tag{8}$$

In the degenerate case it is useful to extend the notion of equation (or solution; we prefer the former): we say that $u$ solves *the integral Poisson equation* if for any $t > 0$, $x \in \mathbb{R}^d$, $u(X_t)$ is $\mathbb{P}_x$-integrable and

$$u(x) = \mathbb{E}_x u(X_t) + \int_0^t \mathbb{E}_x f(X_s)\,ds. \tag{9}$$

This notion is similar to probabilistic or martingale solution of a parabolic equation in [18], also for the degenerate case; in this respect it is worth remembering that a *classical* solution to the degenerate parabolic equation was first constructed by Gikhman [3]. It is also easy to show that a continuous function solution of the integral Poisson equation is a viscosity solution of the Poisson equation, in the sense of [1].

Notice that (9) may be reformulated in the following form: for all $x \in \mathbb{R}^d$,

$$u(X_t) - u(x) + \int_0^t f(X_s)\,ds \text{ is a martingale under } \mathbb{P}_x. \tag{10}$$

Indeed, (10) implies (9) by taking expectation. Vice versa, if we substitute zero by $t'$, and $x$ by $X_{t'}$ in (9) ($t' < t$), then by virtue of the Markov property we get

$$u(X_{t'}) = \mathbb{E}_{X_{t'}} u(X_t) + \mathbb{E}_{X_{t'}} \int_{t'}^t f(X_s)\,ds,$$

or

$$\mathbb{E}\bigg[-u(X_{t'}) + u(X_t) + \int_{t'}^{t} f(X_s)\, ds \bigg| F_{t'}\bigg] = 0.$$

Hence, it follows that

$$\mathbb{E}\bigg[-u(x) + u(X_t) + \int_0^t f(X_s)\, ds \bigg| F_{t'}\bigg] = u(X_{t'}) - u(x) + \int_0^{t'} f(X_s)\, ds,$$

which means exactly the desired martingale property.

Define

$$\tilde{u}(x) = \int_0^\infty |\mathbb{E}_x f(X_t)|\, dt.$$

THEOREM 1. *Let the assumptions $(A_T)$, $(A_b)$ and $(D_\ell)$ be satisfied. We assume that there exists $0 \le \beta$ such that $|f(x)| \le C(1 + |x|^\beta)$ with $C \ge 1$ and that $(A_f)$ holds true. Then (8) defines a continuous function $u$, which is a solution of (9) and satisfies the following properties. For any $m > \beta + 4$, there exists $C_m$ which depends only on $m, \beta$, the value $\sup_{i,x} |b_i(x)|$ and on the constants $C$ in (2), such that*

(11) $$|u(x)| \le \tilde{u}(x) \le C_m(1 + |x|^m), \qquad x \in \mathbb{R}^d,$$

*so that in particular $u$ is $\mu$-integrable. Moreover, again for any $m > \beta + 4$,*

(12) $$\sup_x (1 + |x|^m)^{-1} \bigg| u(x) - \int_0^N \mathbb{E}_x f(X_t)\, dt \bigg| \to 0 \qquad \text{as } N \to \infty.$$

*In addition, $u$ is centered in the sense that*

(13) $$\int u(x) \mu(dx) = 0.$$

*The solution is unique in the class of solutions of (9) which satisfy properties (11) and (13).*

THEOREM 2. *Let the assumptions of Theorem 1 be in force.*

(i) *If there exists $C$ such that*

(14) $$|f(x)| \le C(1 + |x|)^{\beta - 2}$$

*for some $\beta < 0$, then $u$ is bounded. Moreover,*

(15) $$\sup_x |u(x)| \le C \sup_x [|f(x)|(1 + |x|)^{-\beta + 2}],$$

*where the constant $C$ depends only on the constants $C$, $m$, $k$ from (2)–(4) in Proposition 1.*

(ii) *If there exist $C$, $\beta > 0$ with*

(16) $$|f(x)| \leq C(1+|x|)^{\beta-2},$$

*then there exists $C'$ such that*

$$|u(x)| \leq C'(1+|x|)^{\beta}.$$

*Moreover,*

(17) $$\sup_x \frac{|u(x)|}{1+|x|^\beta} \leq C'' \sup_x \frac{|f(x)|}{1+|x|^{\beta-2}},$$

*where the constant $C''$ depends only on the constants $C$, $m$, $k$ from* (2)–(4) *in Proposition* 1.

The assertion of Theorem 1 is used in Theorem 2, which means that the last theorem gives additional information under additional assumptions. Theorem 2 gives in particular a criterion for $u$ to be bounded.

PROOF OF THEOREM 1. The calculations are similar to those in [14]; however, they are not identical. Therefore we present the proof for the reader's convenience.

A. $u$ is well defined and satisfies (11). This follows from [22]; see Proposition 1. Indeed,

$$\tilde{u}(x) = \int_0^\infty |\mathbb{E}_x f(X_t)|\,dt$$
$$= \int_0^\infty \left|\int f(y)\mu_t^x(dy)\right| dt$$
$$= \int_0^\infty \left|\int f(y)[\mu_t^x(dy) - \mu(dy)]\right| dt.$$

Without loss of generality, we assume that $\beta + 2 < m$. Due to the inequalities in Proposition 1, one can choose $p > 1$, $q > 1$ with $p^{-1} + q^{-1} = 1$, such that $p\beta \leq m$ and $q < k+1$.

Indeed, if $\beta = 0$, then it is evident. Consider the case $\beta > 0$. Let $p = m/\beta$. Then $q^{-1} = 1 - \beta/m$, and $(k+1)/q > 1$ is equivalent to $(k+1)(1-\beta/m) > 1$. Since $k+1$ is an arbitrary number less than $m/2$, then the last inequality can be satisfied if $(m/2)(1-\beta/m) > 1$, which is equivalent to $m > \beta + 2$, and this is our assumption. Now, using Hölder's inequality, and denoting all new constants by $C$ (they may be different on each line), one has

$$\int_0^\infty \left|\int f(y)[\mu_t^x(dy) - \mu(dy)]\right| dt$$

$$\leq \int_0^\infty \left(\int |f(y)|^p[\mu_t^x(dy) + \mu(dy)]\right)^{1/p} \left(\int |\mu_t^x - \mu|(dy)\right)^{1/q} dt$$

$$\leq C \int_0^\infty \left(\int (1+|y|^m)[\mu_t^x(dy) + \mu(dy)]\right)^{1/p} (\text{var}(\mu_t^x - \mu))^{1/q} dt$$

$$\leq C \int_0^\infty (1 + \mathbb{E}_x|X_t|^m + \mathbb{E}_\mu|X_t|^m)^{1/p}((1+|x|^m)(1+t)^{-(k+1)})^{1/q} dt$$

$$\leq C(1+|x|^{m'})^{1/q} \int_0^\infty (\mathbb{E}_x|X_t|^m + 1)^{1/p}(1+t)^{-(k+1)/q} dt$$

$$\leq C(1+|x|^{m'})^{1/q} \int_0^\infty (\mathbb{E}_x|X_t|^m)^{1/p}(1+t)^{-(k+1)/q} dt + C(1+|x|^{m'})^{1/q}$$

$$\leq C(1+|x|^{m'}).$$

Thus, $u$ is locally bounded and, moreover, (11) holds true with any $m' > \beta + 4$. The assertion (12) follows from the same calculations with $\int_N^\infty$ instead of $\int_0^\infty$.

B. $u$ satisfies (13). Notice that if some function $g$ is integrable w.r.t. the invariant measure $\mu$, then for any $s > 0$

$$\int \mathbb{E}_x[g(X_s)]\mu(dx) = \int g(x)\mu(dx).$$

Due to (11), the function $\tilde{u}$ is $\mu$-integrable. So, by virtue of Fubini's theorem,

$$\int \int_0^\infty \mathbb{E}_x f(X_s) \, ds \mu(dx) = \int_0^\infty \int \mathbb{E}_x f(X_s) \mu(dx) \, ds.$$

But clearly

$$\int \mathbb{E}_x f(X_s) \mu(dx) = \int f(x)\mu(dx) = 0.$$

C. $u$ is continuous. It follows from the locally uniform convergence (12).

D. $u$ solves the integral Poisson equation (9). Let $t > 0$ be a nonrandom value. First note that

$$u(x) = \int_0^t \mathbb{E}_x f(X_s) \, ds + \int_t^\infty \mathbb{E}_x f(X_s) \, ds,$$

where both integrals are well defined. On the other hand, from the Markov property of $X$,

$$\int_t^\infty \mathbb{E}_x f(X_s) \, ds = \int_0^\infty \mathbb{E}_x \mathbb{E}_{X_t} f(X_s) \, ds$$

$$= \lim_{N \to \infty} \int_0^N \mathbb{E}_x \mathbb{E}_{X_t} f(X_s) \, ds$$

$$= \lim_{N \to \infty} \mathbb{E}_x \int_0^N \mathbb{E}_{X_t} f(X_s) \, ds$$

$$\equiv \lim_{N\to\infty} \mathbb{E}_x u^N(X_t)$$
$$= \mathbb{E}_x u(X_t),$$

where $u^N(x) := \int_0^N \mathbb{E}_x f(X_t)\,dt$. Hence,

$$u(x) - \mathbb{E}_x u(X_t) = \int_0^t \mathbb{E}_x f(X_s)\,ds.$$

This is exactly (9).

E. Uniqueness. For the difference of two solutions, $v = u - u'$, we have due to (9), $v(x) = \mathbb{E}_x v(X_t)$. So

$$v(x) = \mathbb{E}_x v(X_t) \to \int_{\mathbb{R}^d} v(x)\mu(dx) = 0, \qquad t \to \infty.$$

Hence, $v(x) \equiv 0$. $\square$

PROOF OF THEOREM 2. The proof is identical to that in [14]; in particular, the strong Markov property of the process $X_t$ makes possible the use of the formula

$$u(x) = \mathbb{E}_x u(X_{\tau^R}) + \mathbb{E}_x \int_0^{\tau^R} f(X_t)\,dt,$$

which leads to boundedness condition for the function $u$ in the first assertion. We refer to the calculations in [14]. $\square$

**5. Diffusion approximation.** Let $\{X_t, t \geq 0\}$ denote the solution of the SDE

$$dX_t = b(X_t)\,dt + \sigma(X_t)\,dB_t, \qquad X_0 = x,$$

and define $X_t^\varepsilon := X_{t/\varepsilon^2}$, $t \geq 0$. Note that for some Brownian motion $\{B_t^\varepsilon\}$ depending on $\varepsilon$, $X_t^\varepsilon$ solves the SDE

(18) $\qquad dX_t^\varepsilon = \varepsilon^{-2} b(X_t^\varepsilon)\,dt + \varepsilon^{-1} \sigma(X_t^\varepsilon)\,dB_t^\varepsilon, \qquad X_0^\varepsilon = x.$

In this section, we are going to apply Theorem 1 to the singularly perturbed ODE

(19) $\qquad \dfrac{dY_t^\varepsilon}{dt} = F(X_t^\varepsilon, Y_t^\varepsilon) + \varepsilon^{-1} G(X_t^\varepsilon, Y_t^\varepsilon), \qquad 0 \leq t \leq T, Y_0^\varepsilon = y.$

Here $\varepsilon$ is a small parameter. The process $X$ is the same as that of the previous sections, and we will again assume the same conditions $(A_b)$, $(A_T)$ and $(D_\ell)$. $F$ and $G$ are Borel vector-functions. The dimension of $X$ is again $d$, the dimension of $Y$ is $\ell$. We denote again by $L$ the generator of the process $X$. The problem we are interested in is the weak convergence of the slow component $Y^\varepsilon$ as $\varepsilon \to 0$. Concerning (19), we require the Lipschitz condition with respect to the variable $y$, with a constant which may depend on $x$:

$$(A_L) \qquad |F(x,y) - F(x,y')| + |G(x,y) - G(x,y')| \leq C(x)|y - y'|,$$

where $x \to C(x)$ is locally bounded. We now assume that for all $x \in \mathbb{R}^d$, $G(x,\cdot) \in C^1(\mathbb{R}^\ell; \mathbb{R}^\ell)$, that $\partial_y G \in C(\mathbb{R}^{d+\ell}; \mathbb{R}^{\ell^2})$ and the functions $F$, $G$ satisfy the following polynomial growth conditions:

$$(A_P) \qquad \begin{aligned} |F(x,y)| &\leq K(1+|y|)(1+|x|^{q_1}), \\ |G(x,y)| &\leq K(1+|y|)(1+|x|^{q_2}), \\ \|\nabla_y G(x,y)\| &\leq K(1+|x|^{q_3}). \end{aligned}$$

We assume moreover that for all $y \in \mathbb{R}^\ell$ and $j = 1, 2, \ldots, \ell$,

$$(A_G) \qquad \int G_j(x,y)\mu(dx) = 0,$$

where $\mu(dx)$ again denotes the (unique) invariant measure of $X$. It then follows from Theorem 1 that the Poisson equations

$$L\bar{G}_j(x,y) = -G_j(x,y), \qquad j = 1, \ldots, \ell,$$

which we in fact interpret as integral Poisson equations, see (9), have unique centered solutions

$$\bar{G}_j(x,y) = \int_0^\infty E_x G_j(X_t^1, y)\, dt.$$

Moreover, for some $K$ and $q_2'$, $q_3'$, the following holds:

$$(20) \qquad \begin{aligned} |\bar{G}(x,y)| &\leq K(1+|y|)(1+|x|^{q_2'}), \\ \|\nabla_y \bar{G}(x,y)\| &\leq K(1+|x|^{q_3'}). \end{aligned}$$

The values of $q_2'$ and $q_3'$ can be deduced from those of $q_2$ and $q_3$ by using Theorem 1 or Theorem 2, and the fact that $\nabla_y \bar{G} = \overline{\nabla_y G}$.

In the next theorem, we make use of the $S$ topology of Jakubowski on the space $\mathbb{D}([0,T]; \mathbb{R}^\ell)$ of "càdlàg" $\mathbb{R}^\ell$-valued functions defined on $[0,T]$. We refer to [4] for a definition of that topology and the presentation of its properties.

THEOREM 3. *Let the assumptions $(A_b)$, $(D_\ell)$, $(A_T)$ and $(A_L)$, $(A_P)$, $(A_G)$ be satisfied. Then for any $T > 0$, the family of processes $\{Y_t^\varepsilon, 0 \leq t \leq T\}_{0 < \varepsilon \leq 1}$ is uniformly $S$-tight in $\mathbb{D}([0,T]; \mathbb{R}^\ell)$. If $Y$ is an accumulation point of the family $\{Y^\varepsilon, \varepsilon \to 0\}$, then it is a.s. continuous, and it is a solution of the martingale problem associated to the operator*

$$\mathcal{L} = \tfrac{1}{2}\bar{a}_{ij}(y)\partial_{y_i}\partial_{y_j} + \bar{b}_i(y)\partial_{y_i},$$

*where*

$$\bar{b}(y) = \bar{F}(y) + \sum_i \int G_i(x,y) \partial_{y_i} \bar{G}(x,y) \mu(dx)$$

*with*

$$\bar{F}(y) = \int F(x,y) \mu(dx)$$

*and*

$$\bar{a}(y) = \int [G(x,y)\bar{G}^*(x,y) + \bar{G}(x,y)G^*(x,y)] \mu(dx).$$

*If, moreover, the martingale problem associated to $\mathcal{L}$ is well posed (it is easy to state sufficient conditions for that), then $Y^\varepsilon \Rightarrow Y$ in the sense of the S-topology, and $Y$ is the unique (in law) diffusion process with generator $\mathcal{L}$.*

Notice that all integrals in the definition of $\mathcal{L}$ are well defined, as follows from Proposition 1.

PROOF OF THEOREM 3.

*Step* 1. Preliminary computation. Let $f \in C_p^3(\mathbb{R}^\ell)$ (the set of functions of class $C^3$ which, together with their partial derivatives of order 1, 2 and 3, have at most polynomial growth of some order) and define

$$f^\varepsilon(x,y) = f(y) + \varepsilon u(x,y),$$

where $\varepsilon u(x,y)$ is a corrector to $f$, defined as follows. $u$ is the solution of the Poisson equation

$$Lu(x,y) = -\langle \nabla_y f(y), G(x,y) \rangle,$$

or in other words

(21) $$u(x,y) = \langle \nabla_y f(y), \bar{G}(x,y) \rangle,$$

where $\bar{G} : \mathbb{R}^d \times \mathbb{R}^\ell \to \mathbb{R}^\ell$ solves

$$L\bar{G}(x,y) = -G(x,y)$$

in the integral form (9). Note that

$$\int \partial_y G(x,y) \mu(dx) = 0, \qquad y \in \mathbb{R}^\ell,$$

and

$$\partial_y \bar{G}(x,y) = \overline{\partial_y G(x,y)}.$$

For each $\delta > 0$, we associate a mesh $0 = t_0 < t_1 < \cdots < t_n < \cdots$, such that $t_i - t_{i-1} \leq \delta$, $i \geq 0$, and $t_i \to \infty$, as $i \to \infty$. For each $t > 0$, let $N(t)$ denote that smallest integer such that $t \leq t_{N(t)}$. It follows from our definition of the Poisson equation solved by $\bar{G}(x, y)$ that for all $\varepsilon > 0$, $\delta > 0$, the following is a local martingale:

$$M_t^{\varepsilon,\delta} = \sum_{i \leq N(t)-1} \bigg[ \varepsilon u(X_{t_{i+1} \wedge t}^\varepsilon, Y_{t_i}^\varepsilon) - \varepsilon u(X_{t_i}^\varepsilon, Y_{t_i}^\varepsilon)$$

$$+ \frac{1}{\varepsilon} \int_{t_i}^{t_{i+1} \wedge t} \langle \nabla f(Y_{t_i}^\varepsilon), G(X_s^\varepsilon, Y_{t_i}^\varepsilon) \rangle \, ds \bigg].$$

Moreover,

$$\sum_{i \leq N(t)-1} [u(X_{t_{i+1} \wedge t}^\varepsilon, Y_{t_i}^\varepsilon) - u(X_{t_i}^\varepsilon, Y_{t_i}^\varepsilon)]$$

$$= u(X_t^\varepsilon, Y_{t_{N(t)-1}}^\varepsilon) - u(X_0^\varepsilon, Y_0^\varepsilon) - \sum_{i \leq N(t)-2} [u(X_{t_{i+1}}^\varepsilon, Y_{t_{i+1}}^\varepsilon) - u(X_{t_{i+1}}^\varepsilon, Y_{t_i}^\varepsilon)],$$

and for $i \leq N(t) - 2$,

$$\varepsilon u(X_{t_{i+1}}^\varepsilon, Y_{t_{i+1}}^\varepsilon) - \varepsilon u(X_{t_{i+1}}^\varepsilon, Y_{t_i}^\varepsilon)$$

$$= \varepsilon \int_{t_i}^{t_{i+1}} \langle \nabla_y u(X_{t_{i+1}}^\varepsilon, Y_s^\varepsilon), F(X_s^\varepsilon, Y_s^\varepsilon) \rangle \, ds$$

$$+ \int_{t_i}^{t_{i+1}} \langle \nabla_y u(X_{t_{i+1}}^\varepsilon, Y_s^\varepsilon), G(X_s^\varepsilon, Y_s^\varepsilon) \rangle \, ds.$$

Finally,

$$M_t^{\varepsilon,\delta} = \varepsilon u(X_t^\varepsilon, Y_{t_{N(t)-1}}^\varepsilon) - \varepsilon u(X_0^\varepsilon, Y_0^\varepsilon)$$

$$- \varepsilon \sum_{i \leq N(t)-2} \int_{t_i}^{t_{i+1}} \langle \nabla_y u(X_{t_{i+1}}^\varepsilon, Y_s^\varepsilon), F(X_s^\varepsilon, Y_s^\varepsilon) \rangle \, ds$$

$$- \sum_{i \leq N(t)-2} \int_{t_i}^{t_{i+1}} \langle \nabla_y u(X_{t_{i+1}}^\varepsilon, Y_s^\varepsilon), G(X_s^\varepsilon, Y_s^\varepsilon) \rangle \, ds$$

$$+ \frac{1}{\varepsilon} \sum_{i \leq N(t)-1} \int_{t_i}^{t_{i+1} \wedge t} \langle \nabla_y f(Y_{t_i}^\varepsilon), G(X_s^\varepsilon, Y_{t_i}^\varepsilon) \rangle \, ds.$$

We now let $\delta \to 0$ in the last identity, from which we deduce that the following is a local martingale:

$$M_t^\varepsilon = \varepsilon u(X_t^\varepsilon, Y_t^\varepsilon) - \varepsilon u(X_0^\varepsilon, Y_0^\varepsilon)$$

$$- \varepsilon \int_0^t \langle \nabla_y u(X_s^\varepsilon, Y_s^\varepsilon), F(X_s^\varepsilon, Y_s^\varepsilon) \rangle \, ds$$

$$-\int_0^t \langle \nabla_y u(X_s^\varepsilon, Y_s^\varepsilon), G(X_s^\varepsilon, Y_s^\varepsilon)\rangle \, ds$$

$$+\frac{1}{\varepsilon}\int_0^t \langle \nabla_y f(Y_s^\varepsilon), G(X_s^\varepsilon, Y_s^\varepsilon)\rangle \, ds.$$

Moreover, we have that

$$f(Y_t^\varepsilon) = f(Y_0^\varepsilon) + \int_0^t \left\langle \nabla f(Y_s^\varepsilon), F(X_s^\varepsilon, Y_s^\varepsilon) + \frac{1}{\varepsilon} G(X_s^\varepsilon, Y_s^\varepsilon)\right\rangle ds,$$

hence

$$\begin{aligned}
f(Y_t^\varepsilon) &= f(Y_0^\varepsilon) + \int_0^t \langle \nabla f(Y_s^\varepsilon), F(X_s^\varepsilon, Y_s^\varepsilon) + \nabla_y \bar{G}(X_s^\varepsilon, Y_s^\varepsilon) G(X_s^\varepsilon, Y_s^\varepsilon)\rangle \, ds \\
&\quad + \int_0^t \langle \partial^2 f(Y_s^\varepsilon) \bar{G}(X_s^\varepsilon, Y_s^\varepsilon), G(X_s^\varepsilon, Y_s^\varepsilon)\rangle \, ds \\
&\quad + M_t^{\varepsilon, f} \\
&\quad + \varepsilon \langle \nabla_y f(Y_0^\varepsilon), \bar{G}(X_0^\varepsilon, Y_0^\varepsilon)\rangle - \varepsilon \langle \nabla_y f(Y_t^\varepsilon), \bar{G}(X_t^\varepsilon, Y_t^\varepsilon)\rangle \\
&\quad + \varepsilon \int_0^t [\langle \nabla f(Y_s^\varepsilon), \nabla_y \bar{G}(X_s^\varepsilon, Y_s^\varepsilon) F(X_s^\varepsilon, Y_s^\varepsilon)\rangle \\
&\qquad\qquad + \langle \partial^2 f(Y_s^\varepsilon) \bar{G}(X_s^\varepsilon, Y_s^\varepsilon), F(X_s^\varepsilon, Y_s^\varepsilon)\rangle] \, ds,
\end{aligned}$$
(22)

where $\{M_t^{\varepsilon, f}, t \geq 0\}$ is a continuous local martingale which is localized by the sequence of stopping times

$$S_n^\varepsilon := \inf\{t; |Y_t^\varepsilon| > n\}, \qquad n = 1, 2, \ldots.$$

*Step* 2. *S-tightness.* We shall make use of the $S$-topology on $\mathbb{D}([0, T]; \mathbb{R}^\ell)$, introduced by Jakubowski. The following result is a consequence of the results in [4] and [9]:

PROPOSITION 5. *The collection $\{Y_t^\varepsilon, 0 \leq t \leq T\}_{\{0 < \varepsilon \leq 1\}}$ is uniformly S-tight if it satisfies the two conditions*:

(i) *For all $\delta > 0$, there exists $M > 0$ s.t.*

$$P\left(\sup_{0 \leq t \leq T} |Y_t^\varepsilon| > M\right) \leq \delta, \qquad 0 < \varepsilon \leq 1.$$

(ii) $Y_t^\varepsilon - Y_0^\varepsilon = E_t^\varepsilon + V_t^\varepsilon + M_t^\varepsilon$, *with*

(23) $\qquad E_t^\varepsilon \to 0 \qquad$ *in probability, uniformly for $t \in [0, T]$,*

*and for each $n \in \mathbb{N}$,*

(24) $$\sup_{0 < \varepsilon \leq 1} \mathbb{E}(\|V^\varepsilon\|_{T \wedge S_n^\varepsilon} + \langle M^\varepsilon \rangle_{T \wedge S_n^\varepsilon}) < \infty,$$

where $\|V^\varepsilon\|_t$ denotes the total variation of $V^\varepsilon$ between 0 and $t$, and $\langle M^\varepsilon \rangle$ denotes the quadratic variation of the continuous local martingale $M^\varepsilon$.

We first prove that the sequence $(Y_\cdot^\varepsilon)$ satisfies (i). For that sake, we will use (22), with the function $f(y) = \log(1 + |y|^2)$. Recall that the function $u$ depends on $f$. Notice that for this choice of $f$ one has

$$(1 + |y|)|\partial_y f(y)| + (1 + |y|)^2 \|\partial_y^2 f(y)\| + (1 + |y|)^3 \|\partial_y^3 f(y)\| \leq C,$$

and then in particular [see (20)]

$$|u(x,y)| \leq K(1 + |x|^{q_3'}).$$

Consequently, the absolute values of the integrands in the right-hand side of (22) do not exceed $C(1 + |X_s^\varepsilon|^q)$ with some $q < \infty$. So $\{M_t^{f,\varepsilon}\}$ is in fact a martingale, and there exist two constants $C$ and $q$ such that for $0 < \varepsilon \leq 1$,

$$\mathbb{E}\left[\sup_{0 \leq t \leq T} \log(1 + |Y_t^\varepsilon|^2)\right] \leq C \sup_{0 \leq t \leq T} \mathbb{E}(1 + |X_t^\varepsilon|^q) < \infty.$$

This implies that the condition (i) in Proposition 5 is satisfied.

It remains to prove that (ii) is satisfied. For that sake, we choose $f(y) = y$ in (22), yielding

$$Y_t^\varepsilon = Y_0^\varepsilon + E_t^\varepsilon + V_t^\varepsilon + M_t^\varepsilon, \tag{25}$$

where

$$E_t^\varepsilon = \varepsilon \bar{G}(X_0^\varepsilon, Y_0^\varepsilon) - \varepsilon \bar{G}(X_t^\varepsilon, Y_t^\varepsilon),$$

$$V_t^\varepsilon = \int_0^t (I + \varepsilon \nabla_y \bar{G}(X_s^\varepsilon, Y_s^\varepsilon)) F(X_s^\varepsilon, Y_s^\varepsilon) \, ds + \int_0^t \nabla_y \bar{G}(X_s^\varepsilon, Y_s^\varepsilon) G(X_s^\varepsilon, Y_s^\varepsilon) \, ds,$$

and $\{M_t^\varepsilon, t \geq 0\}$ is a continuous local martingale.

Now (23) follows from Corollary 1, (20) and (i), and the first half of (24) follows from $(A_P)$ and (20), and we finally compute $\langle M^\varepsilon \rangle$.

From (22) with $f(y) = |y|^2$,

$$|Y_t^\varepsilon|^2 = |Y_0^\varepsilon|^2 + 2\int_0^t \langle Y_s^\varepsilon, F(X_s^\varepsilon, Y_s^\varepsilon) + \nabla_y \bar{G}(X_s^\varepsilon, Y_s^\varepsilon) G(X_s^\varepsilon, Y_s^\varepsilon) \rangle \, ds$$

$$+ 2\int_0^t \langle \bar{G}(X_s^\varepsilon, Y_s^\varepsilon), G(X_s^\varepsilon, Y_s^\varepsilon) \rangle \, ds$$

$$+ M_t^{\varepsilon,2}$$

$$+ 2\varepsilon \langle Y_0^\varepsilon \bar{G}(X_0^\varepsilon, Y_0^\varepsilon) \rangle - 2\varepsilon \langle Y_t^\varepsilon \bar{G}(X_t^\varepsilon, Y_t^\varepsilon) \rangle$$

$$+ 2\varepsilon \int_0^t [\langle Y_s^\varepsilon, \nabla_y \bar{G}(X_s^\varepsilon, Y_s^\varepsilon) G(X_s^\varepsilon, Y_s^\varepsilon) \rangle + \langle \bar{G}(X_s^\varepsilon, Y_s^\varepsilon), G(X_s^\varepsilon, Y_s^\varepsilon) \rangle] \, ds,$$

where $\{M_t^{\varepsilon,2}, t \geq 0\}$ is a continuous local martingale.

Now from Itô's formula for continuous semimartingales and (25), we deduce that

$$|Y_t^\varepsilon + \varepsilon \bar{G}(X_t^\varepsilon, Y_t^\varepsilon)|^2$$
$$= |Y_0^\varepsilon + \varepsilon \bar{G}(X_0^\varepsilon, Y_0^\varepsilon)|^2$$
$$+ 2\int_0^t \langle Y_s^\varepsilon, F(X_s^\varepsilon, Y_s^\varepsilon) + \nabla_y \bar{G}(X_s^\varepsilon, Y_s^\varepsilon) G(X_s^\varepsilon, Y_s^\varepsilon) \rangle \, ds$$
$$+ 2\varepsilon \int_0^t \langle Y_s^\varepsilon, \nabla_y \bar{G}(X_s^\varepsilon, Y_s^\varepsilon) F(X_s^\varepsilon, Y_s^\varepsilon) \rangle \, ds$$
$$+ 2\varepsilon \int_0^t \langle \bar{G}(X_s^\varepsilon, Y_s^\varepsilon),$$
$$\quad (I + \varepsilon \nabla_y \bar{G}(X_s^\varepsilon, Y_s^\varepsilon)) F(X_s^\varepsilon, Y_s^\varepsilon) + \nabla_y \bar{G}(X_s^\varepsilon, Y_s^\varepsilon) G(X_s^\varepsilon, Y_s^\varepsilon) \rangle \, ds$$
$$+ 2\int_0^t \tilde{Y}_s^\varepsilon \, dM_s^\varepsilon + \langle M^\varepsilon \rangle_t,$$

where $\tilde{Y}_s^\varepsilon = Y_s^\varepsilon + \varepsilon \bar{G}(X_s^\varepsilon, Y_s^\varepsilon)$. Comparing the last two identities, we deduce that

$$\langle M^\varepsilon \rangle_t = 2\int_0^t \langle \bar{G}(X_s^\varepsilon, Y_s^\varepsilon), G(X_s^\varepsilon, Y_s^\varepsilon) \rangle \, ds$$
$$+ \varepsilon^2 |\bar{G}(X_t^\varepsilon, Y_t^\varepsilon)|^2 - \varepsilon^2 |\bar{G}(X_0^\varepsilon, Y_0^\varepsilon)|^2$$
$$+ \varepsilon \int_0^t \psi_\varepsilon(X_s^\varepsilon, Y_s^\varepsilon) \, ds + M_t^{\varepsilon,2} - 2\int_0^t \tilde{Y}_s^\varepsilon \, dM_s^\varepsilon,$$

where

$$|\psi_\varepsilon(x,y)| \leq C(1+\varepsilon)(1+|y|^2)(1+|x|^{3q}).$$

The second half of (24) now follows from (20) and the assumptions on the growth of $G$.

*Step 3.* Identification of the limit. Let $0 \leq s < t \leq T$, and let $\Phi_s$ be a bounded and $S$-continuous functional defined on $\mathbb{D}([0,T]; \mathbb{R}^\ell)$, which is measurable with respect to the $\sigma$-algebra $\sigma(x(r), x \in \mathbb{D}([0,T]; \mathbb{R}^\ell) \ 0 \leq r \leq s)$. Let $f \in C_c^\infty(\mathbb{R}^\ell)$ be a smooth function with compact support. It follows from (22) that for all $a > 0$, such that $t + a < T$,

$$\mathbb{E}([f(Y_{t+a}^\varepsilon) - f(Y_{s+a}^\varepsilon)]\Phi_s(Y))$$
$$= \mathbb{E}\left(\Phi_s(Y) \int_{s+a}^{t+a} \langle \nabla f(Y_r^\varepsilon), F(X_r^\varepsilon, Y_r^\varepsilon) + \nabla_y \bar{G}(X_r^\varepsilon, Y_r^\varepsilon) G(X_r^\varepsilon, Y_r^\varepsilon) \rangle \, dr \right)$$

$$+ \mathbb{E}\bigg(\Phi_s(Y)\int_{s+a}^{t+a}\langle \partial^2 f(Y_r^\varepsilon)\bar{G}(X_r^\varepsilon,Y_r^\varepsilon),G(X_r^\varepsilon,Y_r^\varepsilon)\rangle\,dr\bigg)$$

$$+ \varepsilon\mathbb{E}(\Phi_s(Y)[\langle \nabla_y f(Y_{s+a}^\varepsilon),\bar{G}(X_{s+a}^\varepsilon,Y_{s+a}^\varepsilon)\rangle$$
$$- \langle \nabla_y f(Y_{t+a}^\varepsilon),\bar{G}(X_{t+a}^\varepsilon,Y_{t+a}^\varepsilon)\rangle])$$

$$+ \varepsilon\mathbb{E}\bigg(\Phi_s(Y)\int_{s+a}^{t+a}\langle \nabla f(Y_r^\varepsilon),\nabla_y\bar{G}(X_r^\varepsilon,Y_r^\varepsilon)G(X_r^\varepsilon,Y_r^\varepsilon)\rangle\,dr\bigg)$$

$$+ \varepsilon\mathbb{E}\bigg(\Phi_s(Y)\int_{s+a}^{t+a}\langle \partial^2 f(Y_r^\varepsilon)\bar{G}(X_r^\varepsilon,Y_r^\varepsilon),G(X_r^\varepsilon,Y_r^\varepsilon)\rangle\,dr\bigg).$$

We choose $\delta > 0$ small enough, such that $t + \delta < T$, and deduce from the last identity that

$$\mathbb{E}\bigg(\Phi_s(Y)\int_0^\delta [f(Y_{t+a}^\varepsilon) - f(Y_{s+a}^\varepsilon)]\,da\bigg)$$

$$= \mathbb{E}\bigg(\Phi_s(Y)\int_0^\delta da \int_{s+a}^{t+a}\langle \nabla f(Y_r^\varepsilon),$$

$$F(X_r^\varepsilon,Y_r^\varepsilon) + \nabla_y\bar{G}(X_r^\varepsilon,Y_r^\varepsilon)G(X_r^\varepsilon,Y_r^\varepsilon)\rangle\,dr\bigg)$$

$$+ \mathbb{E}\bigg(\Phi_s(Y)\int_0^\delta da \int_{s+a}^{t+a}\langle \partial^2 f(Y_r^\varepsilon)\bar{G}(X_r^\varepsilon,Y_r^\varepsilon),G(X_r^\varepsilon,Y_r^\varepsilon)\rangle\,dr\bigg)$$

(26)
$$+ \varepsilon\mathbb{E}\bigg(\Phi_s(Y)\int_0^\delta da[\langle \nabla_y f(Y_{s+a}^\varepsilon),\bar{G}(X_{s+a}^\varepsilon,Y_{s+a}^\varepsilon)\rangle$$
$$- \langle \nabla_y f(Y_{t+a}^\varepsilon),\bar{G}(X_{t+a}^\varepsilon,Y_{t+a}^\varepsilon)\rangle]\bigg)$$

$$+ \varepsilon\mathbb{E}\bigg(\Phi_s(Y)\int_0^\delta da \int_{s+a}^{t+a}\langle \nabla f(Y_r^\varepsilon),\nabla_y\bar{G}(X_r^\varepsilon,Y_r^\varepsilon)G(X_r^\varepsilon,Y_r^\varepsilon)\rangle\,dr\bigg)$$

$$+ \varepsilon\mathbb{E}\bigg(\Phi_s(Y)\int_0^\delta da \int_{s+a}^{t+a}\langle \partial^2 f(Y_r^\varepsilon)\bar{G}(X_r^\varepsilon,Y_r^\varepsilon),G(X_r^\varepsilon,Y_r^\varepsilon)\rangle\,dr\bigg).$$

It follows from Lemma 5 in [11] that for any $0 \leq s < t \leq T$,

$$\int_s^t \langle \nabla f(Y_r^\varepsilon), F(X_r^\varepsilon,Y_r^\varepsilon) + \nabla_y\bar{G}(X_r^\varepsilon,Y_r^\varepsilon)G(X_r^\varepsilon,Y_r^\varepsilon) - \bar{b}(Y_r^\varepsilon)\rangle\,dr \to 0$$

and

$$\int_s^t \operatorname{Tr}\partial^2 f(Y_r^\varepsilon)[\bar{G}(X_r^\varepsilon,Y_r^\varepsilon) \otimes G(X_r^\varepsilon,Y_r^\varepsilon) - \tfrac{1}{2}\bar{a}(Y_r^\varepsilon)]\,dr \to 0$$

in probability, as $\varepsilon \to 0$.

We can then take the limit in (26) as $\varepsilon \to 0$, divide by $\delta > 0$, and let $\delta \to 0$ since the process $Y$ is right-continuous, yielding that for all $f \in C_c^\infty(\mathbb{R}^\ell)$, all $0 \le s < t \le T$ and all $\Phi_s$ bounded and $S$-continuous functional defined on $\mathbb{D}([0,T];\mathbb{R}^\ell)$,

$$\mathbb{E}([f(Y_t) - f(Y_s)]\Phi_s(Y)) \tag{27}$$
$$= \mathbb{E}\left(\Phi_s(Y) \int_s^t [\langle \nabla f(Y_r), \bar{b}(Y_r)\rangle + \tfrac{1}{2}\operatorname{Tr} \partial^2 f(Y_r) \bar{a}(Y_r)]\, dr\right),$$

or in other words that

$$M_t^f := f(Y_t) - f(Y_s) - \int_s^t [\langle \nabla f(Y_r), \bar{b}(Y_r)\rangle + \tfrac{1}{2}\operatorname{Tr} \partial^2 f(Y_r) \bar{a}(Y_r)]\, dr$$

is a martingale.

It remains to show that $t \to Y_t$ is a.s. continuous from $[0,T]$ into $\mathbb{R}^\ell$, which is done in the following.

PROPOSITION 6. *Let $\{Y_t, 0 \le t \le T\}$ be an $\ell$-dimensional semimartingale such that for all $1 \le i \le \ell$, all $f \in C_c^\infty(\mathbb{R})$,*

$$M_t^{i,f} := f(Y_t^i) - f(Y_0^i) - \int_0^t [f'(Y_s^i)\bar{b}_i(Y_s) + \tfrac{1}{2} f''(Y_s^i)\bar{a}_{ii}(Y_s)]\, ds$$

*is a martingale. Then $\{Y_t, 0 \le t \le T\}$ is continuous.*

PROOF. We note that the assumption implies that $\forall f \in C^\infty(\mathbb{R})$, $M_t^{i,f}$ is a local martingale. Hence in particular, for each $1 \le i \le \ell$,

$$M_t^i = Y_t^i - Y_0^i - \int_0^t \bar{b}(Y_s^i)\, ds$$

is a local martingale, where $Y_t^i$ denotes the $i$th component of $Y_t$, and it follows from Itô's formula for (possibly discontinuous) semimartingales (see, e.g., [15], page 72) that $\forall f \in C^\infty(\mathbb{R})$,

$$f(Y_t^i) = f(Y_0^i) + \int_0^t f'(Y_s^i)\bar{b}_i(Y_s)\, ds + \int_0^t f'(Y_{s-}^i)\, dM_s^i + \tfrac{1}{2}\int_0^t f''(Y_{s-}^i)\, d[M^i]_s$$
$$+ \sum_{0 < s \le t} (f(Y_s^i) - f(Y_{s-}^i) - f'(Y_{s-}^i)\Delta Y_s^i - \tfrac{1}{2} f''(Y_{s-}^i)(\Delta Y_s^i)^2),$$

where $\{[M^i]_t, t \ge 0\}$ denotes the quadratic variation process of the martingale $M^i$. In the particular case $f(y) = (y^i)^2$, this identity reads

$$(Y_t^i)^2 = (Y_0^i)^2 + 2\int_0^t Y_s^i \bar{b}_i(Y_s)\, ds + 2\int_0^t Y_{s-}^i\, dM_s^i + [M^i]_t.$$

Writing the assumption in the case $f(y) = (y^i)^2$, we obtain that

$$M_t^{i,2} := (Y_t^i)^2 - (Y_0^i)^2 - \int_0^t [2Y_s^i \bar{b}_i(Y_s) + \bar{a}_{ii}(Y_s)]\, ds$$

is a local martingale. Comparing the last two identities, we deduce that $[M^i]_t - \int_0^t \bar{a}_{ii}(Y_s)\,ds$ is a local martingale. Next, comparing the two different ways of writing $(Y_t^i)^3$ and using the identity

$$(Y_s^i)^3 = (Y_{s_-}^i)^3 + 3(Y_{s_-}^i)^2 \Delta Y_s^i + 3Y_{s_-}^i (\Delta Y_s^i)^2 + (\Delta Y_s^i)^3,$$

we deduce that $\sum_{0<s\leq t}(\Delta Y_s^i)^3$ is a local martingale, from which we deduce, by comparing the two different ways of writing $(Y_t^i)^4$ and using the identity

$$(Y_s^i)^4 = (Y_{s_-}^i)^4 + 4(Y_{s_-}^i)^3 \Delta Y_s^i + 6(Y_{s_-}^i)^2 (\Delta Y_s^i)^2 + 4Y_{s_-}^i (\Delta Y_s^i)^3 + (\Delta Y_s^i)^4,$$

that $\sum_{0<s\leq t}(\Delta Y_s^i)^4$ is a local martingale, which is impossible, unless it is identically zero. Since this result holds for any $1 \leq i \leq \ell$, the proposition is established. $\square$

LATP, UMR–CNRS 6632
Centre de Mathématiques
  et d'Informatique
Université de Provence
39 rue F. Joliot Curie
13453 Marseille Cedex 13
France
e-mail: pardoux@cmi.univ-mrs.fr
url: www.latp.univ-mrs.fr/˜pardoux

School of Mathematics
University of Leeds
Woodhouse Lane LS2 9JT, Leeds
United Kingdom
and
Institute of Information Transmission Problems
19 Bolshoy Karetnii 101447
Moscow
Russia
e-mail: veretenn@maths.leeds.ac.uk